\documentclass[11pt]{amsart}
\usepackage{amssymb, amsmath}

\newcommand{\linearmapj}{L_j}
\newcommand{\linearmapi}{L_i}

\newcommand{\linearfamily}{\mathbf{L}}

\newcommand{\R}{\mathbb{R}}

\newcommand{\I}{\mathcal{I}}

\renewcommand{\vec}[1]{\boldsymbol{#1}}
\newcommand{\group}[1]{\mathrm{#1}}

\providecommand{\BL}{\mathop{\rm BL}\nolimits}%

\newcommand{\p}{\mathbf{p}}

\newcommand{\lpar}{\left(}
\newcommand{\rpar}{\right)}

\newcommand{\lip}{\left\langle}
\newcommand{\rip}{\right\rangle}
\newcommand{\lset}{\left\lbrace}
\newcommand{\rset}{\right\rbrace}

\renewcommand{\lpar}{\left(}
\renewcommand{\rpar}{\right)}

\newtheorem{theorem}{Theorem}[section]

\newtheorem{conjecture}[theorem]{Conjecture}
\theoremstyle{remark}
\newtheorem*{remark*}{Remark}

\title{Behaviour of the Brascamp--Lieb constant}
\author{J.~Bennett, N.~Bez,  M.G.~Cowling, and T.C.~Flock}
\thanks{The first and fourth authors were supported by the European Research Council [grant
number 307617], the second author was supported by
JSPS Research Activity Start-up [grant number 26887008] and JSPS Grant-in-Aid for Young Scientists A [grant number 16H05995], and the third author was supported by the Australian Research Council [grant number DP140100531].}

\address{Jonathan Bennett and Taryn C. Flock: School of Mathematics, The Watson Building, University of Birmingham, Edgbaston,
Birmingham, B15 2TT, England.}
\email{J.Bennett@bham.ac.uk, T.C.Flock@bham.ac.uk}
\address{Neal Bez: Department of Mathematics, Graduate School of Science and Engineering,
Saitama University, Saitama 338-8570, Japan}
\email{nealbez@mail.saitama-u.ac.jp}
\address{Michael G. Cowling: School of Mathematics and Statistics, University of New South Wales, UNSW Sydney 2052, Australia}
\email{m.cowling@unsw.edu.au}

\begin{document}
\maketitle
\begin{abstract}
Recent progress in multilinear harmonic analysis naturally raises questions about the local behaviour of the best constant (or bound) in the general Brascamp--Lieb inequality as a function of the underlying linear transformations.
In this paper we prove that this constant is continuous, but is not in general differentiable.
\end{abstract}

%\allowdisplaybreaks

\section{Introduction}
The Brascamp--Lieb inequality is a far-reaching common generalisation of well-known multilinear functional inequalities on euclidean spaces, such as the H\"older, Loomis--Whitney and Young convolution inequalities. It is typically written in the form
\begin{equation}\label{BL}
\int_{H} \prod_{j=1}^m (f_j \circ \linearmapj )^{p_j}
\leq C \prod_{j=1}^m \left(\int_{H_j} f_j\right)^{p_j},
\end{equation}
where $m\in\mathbb{N}$, $H$ and $H_j$ denote euclidean spaces of finite dimensions $n$ and $n_j$ where $n_j\leq n$, equipped with Lebesgue measure for each $1\leq j\leq m$. The functions $f_j:H_j\to\R$ are assumed to be nonnegative.  The maps $\linearmapj:H\to H_j$ are surjective linear transformations, and the exponents $p_j$ satisfy $0\leq p_j\leq 1$.

Following the notation in \cite{BCCT1}
we denote by $\BL(\linearfamily,\p)$
the smallest constant $C$ for which
\eqref{BL} holds for all nonnegative input functions
$f_j\in L^1(H_j)$,
$1\leq j\leq m$.
Here $\linearfamily$
and $\mathbf{p}$ denote
the $m$-tuples $(\linearmapj)_{j=1}^m$
and $(p_j)_{j=1}^m$ respectively. We refer to
$(\linearfamily,\mathbf{p})$
as the \textit{Brascamp--Lieb datum}, and $\BL(\linearfamily,\mathbf{p})$
as the \textit{Brascamp--Lieb constant}.

The Brascamp--Lieb inequality has been studied extensively by many authors, and the delicate questions surrounding the finiteness and attainment of the constant have found some useful answers. Most notably, it was shown by Lieb \cite{Lieb} that in order to compute $\mbox{BL}(\linearfamily,\mathbf{p})$ it is enough to restrict attention to gaussian inputs $f_j$, leading to the expression
\begin{equation}\label{gauss}
\mbox{BL}(\linearfamily,\mathbf{p})=\sup \;\frac{\prod_{j=1}^m(\det A_j)^{p_j/2}}{\det(\sum_{j=1}^m p_j\linearmapj^*A_j\linearmapj)^{1/2}},
\end{equation}
where
%$L^*$ is the transpose of L and
the supremum is taken over all positive definite linear transformations $A_j$ on $H_j$, $1\leq j\leq m$. However, this expression (which of course involves a supremum over a non-compact set) retains considerable complexity in the context of general data. Nevertheless, a concise characterisation of finiteness is available, specifically $
\mbox{BL}(\linearfamily,\mathbf{p})<\infty$ if and only if
\begin{equation}\label{scaling}
 \sum_{j=1}^mp_jn_j=n
 \end{equation}
and
\begin{equation}\label{char}
\dim(V)\leq\sum_{j=1}^mp_j\dim(\linearmapj V)
\end{equation}
hold for all subspaces $V\subseteq H$; see \cite{BCCT1} where a proof of this is given based on \eqref{gauss}, or \cite{BCCT2} for an alternative.

In recent years a variety of generalisations of the Brascamp--Lieb inequality have emerged in harmonic analysis, and have found surprising and diverse applications in areas ranging from combinatorial incidence geometry, to dispersive PDE and number theory -- see for example \cite{BCT,BD}, and perhaps most strikingly \cite{BDG}. These generalisations, which include the multilinear Fourier restriction and Kakeya inequalities (see for example \cite{BBFL} for further discussion), may be viewed as ``perturbations" of the classical Brascamp--Lieb inequality, and naturally raise questions about the local behaviour of $\BL(\linearfamily,\mathbf{p})$ as a function of $\linearfamily$. While seemingly quite innocuous given the formula \eqref{gauss}, very little is known about the regularity of $\BL(\cdot,\mathbf{p})$ in general. In \cite{BBFL} it was shown that the function $\BL(\cdot,\mathbf{p})$ is at least \emph{locally bounded} -- a fact that already has applications in multilinear harmonic analysis and beyond; see \cite{BBFL,BD,BDG}. We note that the analysis in \cite{BBFL} also reveals the relatively simple fact that the finiteness set $\{\linearfamily : \BL(\linearfamily,\mathbf{p})<\infty\}$ is open. The main theorem in this paper is the following:
\begin{theorem}\label{main}
For each $\mathbf{p}$,  the Brascamp--Lieb constant $\BL(\cdot,\mathbf{p})$ is a continuous function.
\end{theorem}
For $\linearfamily_0$ such that $\BL(\linearfamily_0,\mathbf{p})=\infty$, Theorem \ref{main} should be interpreted as ``For any sequence approaching $\linearfamily_0$, the associated Brascamp--Lieb constants approach infinity". Simple, yet instructive examples reveal that this continuity conclusion cannot in general be improved to differentiability.

In contrast, for so-called ``simple" data -- that is, data for which \eqref{char} holds with \emph{strict} inequality for all nontrivial proper subspaces $V$ of $H$ -- Valdimarsson showed in  \cite{Vdiff}, that $\BL(\cdot,\mathbf{p})$ is \emph{differentiable}.

We conclude this section by stating one of the conjectural generalisations of the Brascamp--Lieb inequality that inspired our work. The so-called \emph{nonlinear Brascamp--Lieb inequality} replaces the linear surjections $L_j:\mathbb{R}^n\rightarrow\mathbb{R}^{n_j}$ with \emph{local submersions} $B_j:U\rightarrow\mathbb{R}^{n_j}$, defined on a neighbourhood $U$ of a point $x_0\in\mathbb{R}^n$. In \cite{BB} (see also \cite{BBFL}) it is tentatively conjectured that if $dB_j(x_0)=L_j$ for linear maps $\linearfamily$ such that $\BL(\linearfamily,\mathbf{p})<\infty$, then for some $U'\subseteq U$, there exists $C$ such that 
\begin{equation}\label{nlblconj}
\int_{U'}\prod_{j=1}^m (f_j\circ B_j)^{p_j}\leq C \prod_{j=1}^m\left(\int_{\mathbb{R}^{n_j}}f_j\right)^{p_j}.
\end{equation}
For somewhat restrictive classes of data $(\linearfamily,\mathbf{p})$ this can indeed be achieved (see \cite{BCW,BB,BH} for details and applications), 
while for general data \eqref{nlblconj} is known to hold if the input functions $f_j$ are assumed to have an arbitrarily small amount of regularity (see \cite{BBFL}). Global versions of this inequality are also of interest, see \cite{BBG,KS}.

It is natural to formulate a more quantitative version of this conjecture as follows:
\begin{conjecture}\label{conject}
Given $\varepsilon>0$ there exists $\delta>0$ (depending on the maps $B_j$) such that
%\begin{equation}\label{nlblconjquant}
\[
\int_{B(x_0,\delta)}\prod_{j=1}^m (f_j\circ B_j)^{p_j}\leq (1+\varepsilon)\BL(\linearfamily,\mathbf{p})\prod_{j=1}^m\left(\int_{\mathbb{R}^{n_j}}f_j\right)^{p_j}.
\]
%\end{equation}
\end{conjecture}
This conjecture is intimately related to the continuity of the general Brascamp--Lieb constant. Indeed an elementary scaling and limiting argument (as in \cite{BBG}) shows that, if it were true, then for every $\varepsilon>0$ there would exist $\delta>0$ such that
\begin{equation}\label{imp}
|\BL(\linearfamily(x),\mathbf{p})-\BL(\linearfamily(x_0),\mathbf{p})|<\varepsilon
\end{equation}
whenever $|x-x_0|<\delta$; here $\linearfamily(x)$ denotes the family of linear surjections $(dB_j(x))_j$. (We clarify that one of the two implicit inequalities in \eqref{imp} is a consequence of the elementary \textit{lower} semicontinuity of the Brascamp--Lieb constant.)

One possible application of Conjecture \ref{conject} would be in establishing best constants for local versions of Young's inequality for convolution on noncommutative Lie groups, which is intimately related with best constants for the Hausdorff--Young inequality; this topic has been studied by several authors (see, for example, 
\cite{Beck1,KR,GCMP,CMMP}).
%%Beckner Klein and Russo Garcia-Cuervo, Marco and Parcet, recently by Cowling, Martini, M\"uller and Parcet
 The Baker--Campbell--Hausdorff formula suggests that, for functions supported on small sets, convolution resembles convolution in $\mathbb{R}^n$. If true, Conjecture \ref{conject} would put this on a firm footing.

We refer the reader to \cite{BEsc,BBFL,Zhang} for further discussion and a description of some rather different Kakeya-type and Fourier-analytic generalisations of the Brascamp--Lieb inequality.

Since the writing of this paper we have learnt that Garg, Gurvits, Oliveira and Wigderson have recently discovered a link between Brascamp--Lieb constants and the capacity of quantum operators (as defined in \cite{G04}) which, using ideas from \cite{GGOW}, leads to another proof of the continuity of the Brascamp--Lieb constant in the case of rational data, as well as a quantitative estimate in this case \cite{GGOW2}. We thank Kevin Hughes for bringing this to our attention, and Avi Wigderson for clarifying the connection with this independent work. Lastly, we thank the referee for many thoughtful comments on the original manuscript. 

\subsubsection*{Structure of the paper.\;}

In Section \ref{sec:rank1} we prove Theorem \ref{main} in the relatively straightforward situation where the surjections $L_j$ have rank one (that is, when $n_j=1$ for each $j$) using a well-known formula of Barthe \cite{Barthe}. In Section \ref{sec:genBarthe} we establish a certain general-rank extension of Barthe's formula, which we then combine with the rank-one ideas to conclude Theorem \ref{main} in full generality. Finally, in Section \ref{sec:count}, we present a simple counterexample to the claim that the Brascamp--Lieb constant is everywhere differentiable.

\section{Structure of the proof and the rank-1 case } \label{sec:rank1}

In this section we provide a simple proof of Theorem \ref{main} in the case of rank one maps $L_j$. Our argument is based on the availability of a simpler formula for $\BL(\linearfamily,\mathbf{p})$ due to Barthe \cite{Barthe} in that case. This is a natural starting point since our proof in the general-rank case proceeds by first establishing a suitable Barthe-type formula which holds in full generality.

\begin{proof} We begin by setting up some notation. As each $\linearmapj$ has rank one, there exists a vector $v_j$ such that $\linearmapj(x) = \langle v_j,x\rangle $. We denote by $\mathbf{v}$ the $m$-tuple $(v_j)_{j=1}^m$, and identify $\mathbf{v}$ with $\linearfamily$. Set  \[ \I=\{I: I\subseteq\{1,2,...,m\}, |I|=n\}.\]
For each $I\in\I$, set $p_I= \prod_{i\in I} p_i$ and
\[ d_I=\det((v_i)_{i\in I})^2,\]
where the determinant of a sequence of $n$ vectors in $\R^n$ is the determinant of the $n\times n$ matrix whose $i$th column is the $i$th term in the sequence. Finally, let $\vec{d}=(d_I)_{\I}$.  We view $\vec{d}$ as an element in $\R^N$ where $N={{n}\choose{m}}$. 

Barthe's formula \cite{Barthe} for the best constant is
\begin{equation}\label{rhs} \BL(\mathbf{v},\p)^2 = \sup_{\lambda_i>0} \frac{ \prod_{i=1}^m \lambda_i^{p_i} }{ \sum_{\mathcal{I}} d_I p_I \lambda_I},
\end{equation}
where  $\lambda_I= \prod_{i\in I} \lambda_i$.  By definition, for each $I$, $d_I$ is a continuous function of $\linearfamily$, and so it is enough to show that the right-hand side of \eqref{rhs} is a continuous function of $\vec{d}$.   

First, lower semicontinuity is immediate from the definition, as a supremum of lower semicontinuous functions is itself lower semicontinuous. It thus suffices to prove upper semicontinuity. 

Fix a point $\widetilde{\vec{d}}\in\R^N$. If for each $I$, $\widetilde{d}_I=0$, then the right-hand side of \eqref{rhs} is infinite and uppersemicontinuity is immediate.  If not, then set $D=\min_{I:\widetilde{d}_I\neq 0} \widetilde{d}_I>0$.  For $\delta\in(0, D)$ and $\vec{d}\in\R^N$ such that $|\vec{d}-\widetilde{\vec{d}}|\leq\delta$,
\[ \sup_{\lambda_i>0} \frac{ \prod_{i=1}^m \lambda_i^{p_i} }{ \sum_{\mathcal{I}} d_I p_I\lambda_I} \leq  \sup_{\lambda_i>0} \frac{ \prod_{i=1}^m \lambda_i^{p_i} }{ \sum_{\mathcal{I}: \widetilde{d}_I\neq0} d_I p_I\lambda_I} \leq \sup_{\lambda_i>0} \frac{ \prod_{i=1}^m \lambda_i^{p_i} }{ \sum_{\mathcal{I}: \widetilde{d}_I\neq0} (\widetilde{d}_I-\delta) p_I\lambda_I},\]
and so
\[ \sup_{\lambda_i>0} \frac{ \prod_{i=1}^m \lambda_i^{p_i} }{ \sum_{\mathcal{I}} d_I p_I\lambda_I} \leq \sup_{\lambda_i>0} \frac{ \prod_{i=1}^m \lambda_i^{p_i} }{ \sum_{\mathcal{I}: \widetilde{d}_I\neq0}\widetilde{d}_Ip_I\lambda_I}\left( 1- \delta \frac{\sum_{\mathcal{I}:\widetilde{d}_I\neq0}p_I\lambda_I}{\sum_{\mathcal{I}:\widetilde{d}_I\neq0} \widetilde{d}_Ip_I\lambda_I}  \right)^{-1}.\]
Focusing on the second term in the product, for each $I$ such that $\widetilde{d}_I\neq 0$,
\[ \frac{p_I\lambda_I}{\widetilde{d}_Ip_I\lambda_I} \leq \frac{1}{D}.\]   Thus the second term is bounded by $ (1-\delta/D)^{-1}$, and so
\[ \lim_{\delta\to0} \sup_{\lambda_i>0} \frac{ \prod_{i=1}^m \lambda_i^{p_i} }{ \sum_{\mathcal{I}} d_I p_I\lambda_I} \leq  \sup_{\lambda_i>0} \frac{ \prod_{i=1}^m \lambda_i^{p_i} }{ \sum_{\mathcal{I}: \widetilde{d}_I\neq0}\widetilde{d}_Ip_I\lambda_I}\leq\sup_{\lambda_i>0} \frac{ \prod_{i=1}^m \lambda_i^{p_i} }{ \sum_{\mathcal{I}}\widetilde{d}_Ip_I\lambda_I}.\]
This proves the required upper semicontinuity.
\end{proof}

\section{A generalisation of Barthe's formula and the proof of Theorem \ref{main}}\label{sec:genBarthe}

To extend the proof from Section \ref{sec:rank1} to the general case we first need an analogue of Barthe's formula \eqref{rhs} for the best constant \eqref{gauss}. 
A key step in obtaining such a formula is the parametrisation of positive definite matrices by a rotation matrix and a diagonal matrix of their (positive) eigenvalues. In the rank one case this is equivalent to Barthe's parametrisation. In \cite{Vdiff} Valdimarsson uses a related approach, although he parametrises positive definite matrices by symmetric matrices, rather than by diagonal matrices and rotations.

It will be helpful to simplify notation so as to avoid double sums of the form $\sum_{i=1}^m \sum_{j = 1}^{n_i} a_{ij}$.
To do so, we define $K = \sum_{i=1}^m n_i$, and write $a_{ij}$ as $a_{k}$, where 
$k=n_0+\cdots +n_{i-1}+j$, $n_0=0$, and $1\leq k\leq K$.
%$i$ is such that
%\[
%n_1 + \dots + n_{i-1} < k \leq n_1 + \dots + n_i,
%\]
%and $j = k - n_1 + \dots + n_{i-1}$.
Given this relationship between $(i,j)$ and $k$, define $q_k=p_i$, so that $(q_1,\hdots,q_K)$ is a $K$-tuple whose first $n_1$ entries are $p_1$, next $n_2$ entries are $p_2$, and so on.

We write $\I$ for the set of all subsets $I$ of $\lset 1,  \dots, K\rset$ of cardinality $n$. For each $I\in\I$ define $q_I=\prod_{k\in I} q_k$.  Similarly, given a family $\lambda_1,\hdots,\lambda_K>0$ we define $\lambda_I=\prod_{k\in I} \lambda_k$. In what follows $R_i$ will denote a rotation on $\mathbb{R}^{n_i}$ for each $1\leq i\leq m$, and we denote by $\mathbf{R}$ the $m$-tuple $(R_i)_{i=1}^m$.
\begin{theorem} \label{thm:MGCformula}%Given a Brascamp-Lieb datum $(\linearfamily, \mathbf{p})$, the Brascamp--Lieb constant $\BL(\linearfamily, \mathbf{p})$ is given by
\[
\BL(\linearfamily,\mathbf{p})^{2} = \sup\lset
\frac{\prod_{k=1}^K \lambda_k^{q_k} }{ \sum_{I\in \I}  \lambda_I q_I d_I }   : \lambda_k \in \R^+, R_i \in \group{SO}(n_i)
\rset,
\]
where $d_I=d_I(\linearfamily,\mathbf{R})$ is a nonnegative continuous function for each $I\in\I$.
\end{theorem}
\begin{remark*}
The functions $d_I(\linearfamily,\mathbf{R})$ will be specified in the proof. In the rank-1 case,  each $R_i=1$ and $d_I=\det((v_i)_{i\in I})^2$ as before.
\end{remark*}
\begin{proof}
Recall from \eqref{gauss},
\[
\mbox{BL}(\linearfamily,\mathbf{p})=\sup \;\frac{\prod_{i=1}^m(\det A_i)^{p_i/2}}{\det\left(\sum_{i=1}^m p_i\linearmapi^*A_i\linearmapi\right)^{1/2}}.\]
Here $A_i$ is a positive definite $n_i \times n_i$ matrix. Further, $A_i$ is of the form $R_i^* D_i R_i$, where $R_i$ is a rotation matrix, with transpose $R_i^*$, and $D_i$ is a diagonal matrix, with positive diagonal entries $\lambda_i^1, \dots, \lambda_i^{n_i}$. Using the notation introduced above,
\[ \prod_{i=1}^m(\det A_i)^{p_i/2} = \prod_{k=1}^K \lambda_k^{q_k/2} . \]

Let the column vectors $\{ \vec{e}_i^j : j = 1, \dots, n_i \}$ be the standard basis for $\R^{n_i}$, so that $D_i = \sum_{j=1}^{n_i} \lambda_i^j \vec{e}_i^j({\vec{e}_i^j})^*$, and let $v_i^j = \linearmapi^* R_i^* \vec{e}_i^j$. Then,
\begin{align*}
\sum_{i=1}^m p_i\linearmapi^*A_i\linearmapi & =  \sum_{i=1}^m p_i\linearmapi^*R_i^* D_i R_i\linearmapi \\
& =  \sum_{i=1}^m p_i\linearmapi^*R_i^*\left(  \sum_{j=1}^{n_i} \lambda_i^j \vec{e}_i^j(\vec{e}_i^j)^* \right)   R_i\linearmapi \\
& =  \sum_{i=1}^m p_i \sum_{j=1}^{n_i} \lambda_i^j \linearmapi^* R_i^* \vec{e}_i^j \left(\linearmapi^* R_i^* \vec{e}_i^j\right)^* \\
& =  \sum_{k=1}^K q_k \lambda_k v_k  v_k^*
\end{align*}
where the $v_k$ are the $v_i^j$ in our chosen order, the $\lambda_k$ are the corresponding $\lambda_i^j$, and the $q_k$ are the corresponding $p_i$ (as described above). Set
\[ T = \sum_{k=1}^K q_k \lambda_k v_k  v_k^* .\]
To compute $\det(T)$, we follow Barthe, and use the  Cauchy--Binet formula.
Define the $n \times K$ matrices $A$ and $B$ to be the matrices whose $k$th columns are the vectors $\lambda_k  q_k v_k$ and $v_k$, respectively and the $K \times n$ matrix $C$ to be $B^*$, the matrix whose $k$th row is the vector $v_k^*$. Recall that $\I$ denotes the set of all subsets of $\{1, \dots, K\}$ of cardinality $n$; write $A_I$ and $B_I$ for the $n \times n$ matrices whose columns are the vectors $\lambda_k  q_k v_k$ and $v_k$, respectively where $k \in I$, and $C_I=B_I^*$. % for the $n \times n$ matrix whose rows are the vectors $v_k^*$ where $k \in I$.
Then
\[
\det(T)=\det(AC)
%= \det\lpar \sum_{k=1}^K \lambda_k v_k \otimes v_k^* \rpar
= \sum_{I\in \I} \det (A_I C_I)
=  \sum_{I\in \I} \lpar \prod_{k \in I} \lambda_k  q_k \rpar \det ({B_I} C_I)
= \sum_{I\in \I}  \lambda_I q_I d_I.
\]
Here $\lambda_I = \prod_{k \in I} \lambda_k$, $q_I = \prod_{k \in I} q_k$ and $d_I = \det (B_I C_I)$. Evaluating $d_I$ using the definition of $B_I$ and $C_I$ yields $ d_I= \det((v_k)_{k\in I})^2$ where $v_k=v_{i}^j= \linearmapi^* R_i^* \vec{e}_i^j$.

We conclude that
% \begin{align*}
% \BL(\linearfamily,\mathbf{p})
% &= \sup\lset \frac{ \int_{\R^n} \prod_{i=1}^m \lpar\gamma_i \circ \linearmapi (x)\rpar^{p_i}  \,dx }{ \prod_{i=1}^m \lpar \int_{\R^{n_i}} \gamma_i(x) \,dx \rpar^{p_i} } :
% \text{$\gamma_i$ are gaussians}  \rset \\
% &= \sup\lset \frac{ \det(T)^{-1/2} }{ \prod_{k=1}^K \lambda_k^{-q_k/2} }  : \lambda_k \in \R^+, R_i \in \group{SO}(n_i) \rset  \\
% &=\left( \sup\lset
% \frac{\prod_{k=1}^K \lambda_k^{q_k} }{ \sum_{I\in \I}  \lambda_I q_I d_I }   : \lambda_k \in \R^+, R_i \in \group{SO}(n_i)
% \rset \right)^{1/2}
% \end{align*}
% and hence
\begin{equation*}%\label{MGCformula}
\BL(\linearfamily,\mathbf{p})^{2} = \sup\lset
\frac{\prod_{k=1}^K \lambda_k^{q_k} }{ \sum_{I\in \I}  \lambda_I q_I d_I }   : \lambda_k \in \R^+, R_i \in \group{SO}(n_i)
\rset,
\end{equation*}
where each $d_I$ is manifestly nonnegative and continuous as a function of $\linearfamily$ and $\mathbf{R}$.
\end{proof}
%\section{Proof of Theorem \ref{main}}\label{sec:shortcut}
The argument presented in the rank-1 case in Section \ref{sec:rank1}, combined with Theorem \ref{thm:MGCformula}, quickly leads to Theorem \ref{main}. Define the function $F(\vec{d})$ by
\[
F(\vec{d})=\sup_{\lambda_k>0}\frac{ \prod_{k = 1}^K \lambda_k^{q_k}  }{ \sum_{I \in \I} d_I q_I \lambda_I }.
\]
Then
\[
\BL(\linearfamily,\mathbf{p})^{2} = \sup\lset
F(\vec{d}(\linearfamily,\mathbf{R})) :  R_i \in \group{SO}(n_i)
\rset.
\]
Now, $F(\vec{d})$ is continuous by the rank-1 argument in Section \ref{sec:rank1}.  As the supremum
(as the parameter varies) of a family of continuous functions
continuously parametrised by a compact set is continuous, the Brascamp--Lieb constant is as well.

\section{A non-differentiable example}\label{sec:count}
In general, Brascamp--Lieb constants are surprisingly difficult to compute.  However, consider the general 4-linear rank-1 case when $\mathbf{p}=(1/2,1/2,1/2,1/2)$:
\[
\int \prod_{i=1}^4 f_i^{1/2}(\lip x,v_i\rip)\;dx\leq C  \prod_{i=1}^4 \left(\int f_i \right)^{1/2}.
\]
Here we can exploit symmetry to compute:
\begin{equation}\label{gennondiff}
\mbox{BL}(\mathbf{v}, \mathbf{p})^2 =  2\left({|\det(v_1v_2)\det(v_3v_4)| + |\det(v_1v_3)\det(v_2v_4)|+|\det(v_1v_4)\det(v_2v_3)| }\right)^{-1},
\end{equation}
where again the determinant of a collection of $2$ vectors in $\R^2$ is the determinant of the $2\times2$ matrix with columns $v_i$.

Indeed, by a change of variables and a rescaling, it is enough to consider the following case:
 \[
 \int\int f_1^{1/2}(x)f_2^{1/2}(y)f_3^{1/2}(x-y)f_4^{1/2}(x+ay) \;dxdy \leq C  \prod_{i=1}^4 \left(\int f_i \right)^{1/2}.
\]
Using the Cauchy--Schwarz inequality and a change of variables,
 \begin{multline*}
 \int\int f_1^{1/2}(x)f_2^{1/2}(y)f_3^{1/2}(x-y)f_4^{1/2}(x+ay)\;dxdy\\
 \leq  \left(\int\int f_1(x)f_2(y)\;dxdy\right)^{1/2}\left( \int\int f_3(x-y)f_4(x+ay)\;dxdy\right)^{1/2}\\
  \leq |a+1|^{-1/2} \prod_{i=1}^4 \left(  \int f_i \right)^{1/2}.
 \end{multline*}
This inequality is sharp whenever the application of the Cauchy--Schwarz inequality is sharp, in this case when $a>0$, as may be seen by considering suitable gaussians. Repeating this argument for each possible pairing of functions yields the following formula for the sharp constant when $a\neq 0,-1$:
\[ \min\{ 1, |a|^{-1/2}, |a+1|^{-1/2} \}\]
or, equivalently, 
  \[\left( \frac{2}{|a|+|a+1|+1}\right)^\frac{1}{2} .\]
By Theorem \ref{main} this formula must hold for all $a$. %(Alternately, one can compute the Brascamp--Lieb constant in these cases by factoring through the critical subspace.)
Changing variables back to the general setting yields the general formula \eqref{gennondiff}.
\begin{remark*}
We suspect that in general the dependence of $\BL(\linearfamily,\p)$ on $\linearfamily$ is at least locally H\"older continuous,
and we hope to return to this in a subsequent paper.
\end{remark*}

%	\affiliationone{Jonathan Bennett and Taryn C. Flock \\
%	 School of Mathematics\\
%	%The Watson Building\\
%	University of Birmingham\\
%	 Edgbaston\\
%	Birmingham\\
%	 B15 2TT\\
%	 England
%	\email{J.Bennett@bham.ac.uk\\ T.C.Flock@bham.ac.uk} }
%	\affiliationtwo{ Neal Bez \\
%	Department of Mathematics,\\
%	Graduate School of Science and Engineering\\
%	Saitama University\\
%	Saitama\\
%	 338-8570 \\
%	 Japan
%	\email{nealbez@mail.saitama-u.ac.jp} }
%	\affiliationthree{Michael G. Cowling\\
%	School of Mathematics and Statistics\\
%	University of New South Wales\\
%	Sydney \\
%	2052\\
%	 Australia
%	\email{m.cowling@unsw.edu.au} }

\end{document}